\documentclass[12pt]{amsart}
\usepackage{templatestyle}
\usepackage{templatecommand}
\usepackage{templatetensor}
\usepackage{comment}

\title[Negative holomorphic sectional curvature and the K\"ahler-Ricci flow]%top of the pages, abbreviate if neccesary
{K\"ahler manifolds with negative holomorphic sectional curvature, K\"ahler-Ricci flow approach} %title

\author%[Ryosuke Nomura] %top of the pages
{Ryosuke Nomura}%full name

\address{Graduate School of Mathematical Sciences, The University of Tokyo \endgraf
	3-8-1 Komaba, Meguro-ku, Tokyo, 153-8914, Japan.}

\email{nomu@ms.u-tokyo.ac.jp}

\thanks{Classification AMS 2010: 53C55, %Differential geometry/Hermitian and Kahlerian manifolds
	32W20.%Several complex variables and analytic spaces/Complex Monge-Ampere operators
}


\begin{document}
\begin{abstract}
	Recently, Wu-Yau and Tosatti-Yang established the connection between  the negativity of  holomorphic sectional curvatures and the positivity of canonical bundles for  compact K\"ahler manifolds. In this short note, we give another proof of their theorems by using the K\"ahler-Ricci flow.
\end{abstract}

\maketitle

\section{Introduction}

In this note, we provide a \krf  \ approach to the following two theorems, which represent the relationship between the negativity of the holomorphic sectional curvature and the positivity of the canonical bundle $K_X$ of a compact \kahler \ manifold $X$.

\begin{theorem}[{\cite[Theorem 2]{WuYau16}, \cite[Corollary 1.3]{TosattiYang15}}]\label{amplethm}
If $X$ admits a \kahler \ form with strictly negative holomorphic sectional curvature, then the canonical bundle $K_X$ is ample. In particular, $X$ is projective.
\end{theorem}

\begin{theorem}[{\cite[Theorem 1.1]{TosattiYang15}}]\label{nefthm}
If $X$ admits a \kahler \ form with semi-negative holomorphic sectional curvature, then the canonical bundle $K_X$ is nef.
\end{theorem}

The original proofs of both theorems are based on the following idea, in \cite{WuYau16}, constructing a \kahler \ form $\ome \in 2\pi c_1(K_X) + \varepsilon [\omh]$ satisfying  
\begin{align*}
\ric(\ome ) = - \ome + \varepsilon \omh,
\end{align*}
and considering the limiting behavior of $\ome $ as $\varepsilon \searrow 0$. Here, $\omh $ is a \kahler \ form whose holomorphic sectional curvature is (strictly/semi-) negative.
The objective of this note is to simplify the proofs by replacing $\omega_\varepsilon$ by the \krf \ $\omt$. 

These theorems are originated from the conjecture of Yau (see \cite[Conjecture 1.2]{HeierLuWong16}).
For a historical background, we refer to \cite{HeierLuWong16,WuYau16,WuYau16remark,TosattiYang15,2016arXiv160601381D} and the references therein. 

We remark that  Diverio and Trapani \cite{2016arXiv160601381D} showed that the ampleness of $K_X$ can be obtained under the assumption that the holomorphic sectional curvature is semi-negative everywhere and strictly negative at one point.  
For the moment, we can only prove the above two theorems.

\noindent 
\textbf{Acknowledgment. }The author would like to thank his supervisor Prof. Shigeharu Takayama for various comments.
This work is supported by the Program for Leading Graduate Schools, MEXT, Japan.

\section{Properties of the \krf }

In this section, we summarize the basic properties of the (normalized) \krf \ which will be used later. For more detailed exposition, we refer to the book \cite{BEG13IntrotoKRF}.
In the following argument, we will denote by $X$ a compact \kahler \  manifold of dimension $n$. 

\begin{definition}
A smooth family of \kahler \ forms $\{ \omt \}_{t\ge 0}$ is called
the \textit{normalized \krf} starting from $\omz$ 
 if it satisfies the following equation:
\begin{align}\label{krfeq}
\begin{cases}
	\dt \omt
		&= -\ric(\omt) - \omt, \\[5pt]
	\omt|_{t=0}
		&=\omz.
\end{cases}
\end{align}
\end{definition}

By considering the cohomology class in $H^{1,1}(X, \rr)$ of (\refs{krfeq}), 
the normalized \krf \ $\omt$ belongs to the cohomology class $\alpha_t \in H^{1,1 }(X, \rr)$ defined by
\begin{align}\label{cohom}
	\alpha_t  \deq 
			e^{-t }[\omz ] + (1-e^{-t}) 2\pi  c_1(K_X) .
\end{align}
The  maximal existence theorem for the \krf \ is stated as follows. 
\begin{theorem}[{\cite{TianZhang06KRFProjGenType,Zhang06thesis}}, see also {\cite[3.3.1]{BEG13IntrotoKRF}}]\label{existthm}
For any \kahler \ form $\omz$,  
the normalized \krf
\ $\omt$ starting from $\omz$ exists uniquely for $t \in [0,T)$ and cannot extend beyond $T$.  Here, $T$ is defined by
\begin{align}\label{maximaltime}
	T \deq \sup \{ t>0 \mid \alpha_t \mbox{ defined by (\refs{cohom}) is a \kahler \ class}  \},
\end{align}
and called the maximal existence time.
In particular, $\omt$ exists for $t\in [0,\infty)$ if and only if $K_X$ is nef, i.e. $2\pi c_1(K_X)$ belongs to the closure of the \kahler \ cone of $X$.
\end{theorem}

We need the parabolic Schwarz lemma obtained by Song-Tian \cite{SongTian07KRFSurfKod} applied to the identity map (see also \cite[3.2.6]{BEG13IntrotoKRF}). This is a parabolic analogue of the Schwarz lemma due to Yau \cite{Yau78Schwarz}. 
\begin{proposition}\label{schwarzprop}
Let $\omt$ be 
the normalized \krf
\ and $\omh$ be an arbitrary \kahler \ form.
Then we have the following inequality:
\begin{align*}
\dalt \log \tromt(\omh ) 
	\le 1
		+ \dfrac{ \guijb (t) \guklb (t) ?{\widehat{R}}_i\jbar k \lbar ? }{ \tromt (\omh ) }.
\end{align*}
\end{proposition}

The next Proposition due to Royden \cite[Lemma]{Royden80} (see also \cite[Lemma 2.1]{WongWuYau12}) will 
relate the negativity assumption on the holomorphic sectional curvature with
a uniform estimate for the normalized \krf \ via the parabolic Schwarz lemma.
\begin{proposition}\label{Roydenprop}
Let $\omh$ be a \kahler \ form on $X$, and denote by $\widehat{H} $ the holomorphic sectional curvature of $\omh$. Assume that there exists a non-negative constant $\kappa \ge 0$ such that for any tangent vector $\xi \in T^{1,0}X$, we have
\begin{align}\label{holsect}
\widehat{H} (\xi ) \le -\kappa \le 0.
\end{align}
Then, for any \kahler \ form $\omega$, we have
\begin{align*}
\guijb \guklb ?{\widehat{R}}_i\jbar k \lbar ? \le -\kappa \dfrac{n+1}{2n} \left( \trom(\omh )\right)^2 \le 0,
\end{align*} 
where $\omega = \ii \gijb \dzidzjb $ and $?{\widehat{R}}_i\jbar k \lbar ? $ is the curvature tensor of $\omh$.
\end{proposition}

\section{Proof of Theorems  via \krf % \refs{amplethm} and Theorem \refs{nefthm} 
}
\begin{proof}[Proof of Theorem \refs{nefthm}]
Let $\omh$ be a \kahler \ form whose holomorphic sectional curvature is semi-negative i.e. $\kappa= 0$ in (\refs{holsect}). 
Let $\omt $ be the normalized \krf \ starting from arbitrary \kahler \ form $\omz$ on $X$.  
By Theorem \refs{existthm}, the nefness of $K_X$ is equivalent to the long time existence of $\omt$. By definition of the maximal existence time (\refs{maximaltime}) and Theorem \refs{existthm}, it  is enough to show that if $\omt$ exists for $[0,T_0)$ with $T_0<\infty$, then $\alpha_{T_0}$ is  a \kahler \  class. 

By Proposition \refs{schwarzprop} and, Proposition \refs{Roydenprop} we have 
\begin{align*}
\dalt \left( \log \tromt(\omh ) - t \right)
	\le \dfrac{ \guijb (t) \guklb (t) ?{\widehat{R}}_i\jbar k \lbar ? }{ \tromt (\omh ) } \le 0.
\end{align*}
Applying the maximum principle, for any $t \in [0,T_0)$, we have 
\begin{align*}
	\tromt(\omh) 
		\le e^t  \max_X \tromz(\omh )
		\le e^{T_0} \max_X \tromz(\omh ) \ddeq C.
\end{align*}
This is equivalent to the following inequality which holds for any $t \in [0,T_0)$:
\begin{align}\label{C2ineq}
	 \dfrac{1}{C}\omh \le \omt. 
\end{align}
Therefore, for any irreducible subvariety $V\subset X$ of positive dimension, the intersection number can be estimated as follows:
\begin{align*}
	\int_V [\alpha_{T_0}]^{\dim V} = \lim_{t\nearrow T_0 }\int_V  \omt ^{\dim V}\ge \dfrac{1}{C^{\dim V}}\int_V \omh^{\dim V}>0.
\end{align*}
By Demailly-P\u aun's characterization of the \kahler \ cone \cite[Main Theorem 0.1]{DemaillyPaun04}, the limiting class $\alpha_{T_0}$ is  \kahler. 
\end{proof}

The idea of avoiding higher order estimates by using the Demailly-P\u aun's theorem can be found in the proof of \cite[Theorem 1.1]{Zhang10ScalFTSingKRF}.
\begin{remark}\label{higherrem}
We can also prove that $\omt$ converges to a smooth \kahler \ form as $t\nearrow T_0$, in particular $\alpha_{T_0}$ is  \kahler.
In fact, since the volume form $\omt^n$ is uniformly bounded from above (see \cite[3.2.3]{BEG13IntrotoKRF}),
(\refs{C2ineq}) gives a uniform estimate for $\omt$:
\begin{align}\label{C2est}
	 \dfrac{1}{C}\omh \le \omt \le C^\prime \omh.
\end{align}
Therefore we obtain the higher order estimates (see for example \cite[3.2.16]{BEG13IntrotoKRF}), which guarantees the convergence.
\end{remark}

\begin{proof}[Proof of Theorem \refs{amplethm}]
%By the assumption in Theorem \refs{amplethm}, there exists  a \kahler \ form $\omh$ whose holomorphic sectional curvature is strictly negative i.e. $\kappa>0$ in (\refs{holsect}). Let $\omt $ be the normalized \krf \ starting from arbitrary \kahler \ form $\omz$ on $X$.
Let $\omh$ be a \kahler \ form whose holomorphic sectional curvature is  strictly negative i.e. $\kappa>0$ in (\refs{holsect}). 
By Theorem \refs{nefthm}, $K_X$ is nef, and therefore $\omt$ exists for $t \in [0,\infty)$.
%It is enough to show that there exists a sequence $\{t_i \}\subset [0,\infty)$ such that $t_i $ tends to $\infty$ and $\omti$ converges to some \kahler \ form $\omega_\infty $ as $i \rightarrow \infty$. 
%In fact, since $\alpha_t = [\omt ]$ converges to $2\pi c_1(K_X)$ as $t \rightarrow \infty$, the \kahler \ form $\omega_\infty$ represents $2\pi c_1(K_X)$, and therefore $K_X$ is ample. 

By  Proposition \refs{schwarzprop} and Proposition \refs{Roydenprop}, we get 
\begin{align*}
\dalt \log \tromt(\omh ) 
	\le 1+\dfrac{ \guijb (t) \guklb (t) ?{\widehat{R}}_i\jbar k \lbar ? }{ \tromt (\omh ) }
	\le 1-\kappa \dfrac{n+1}{2n}\tromt(\omh ) .
\end{align*}
Applying the maximum principle, for any $t \in [0,\infty )$, we have $\tromt(\omh) \le C$ where
\begin{align*}
 C \deq \max \left\{\dfrac{2n}{\kappa (n+1)}, \max_X \tromz(\omh)  \right\}>0.
\end{align*}
This  gives, for any $t \in [0,\infty )$,   
\begin{align}\label{C2ineq2}
	 \dfrac{1}{C}\omh \le \omt.
\end{align}
Since $\alpha_t = [\omt ] $ converges to $2\pi c_1(K_X)$ as $t \rightarrow \infty$, the same argument as in the proof of Theorem \refs{nefthm} shows the ampleness of $K_X$.
\end{proof}

\begin{remark}
Tian-Zhang \cite{TianZhang06KRFProjGenType}, and Zhang \cite{Zhang06thesis} showed that under the assumption on the ampleness of $K_X$, any normalized \krf \ $\omega_t$ converges to the \ke \ metric with negative Ricci curvature. 
\end{remark}

\bibliographystyle{amsalphaurlmod}
\bibliography{reference}
\end{document}